\documentclass[12pt]{amsart}

\usepackage{amsfonts}
\usepackage{amscd}
\usepackage{amssymb}

\allowdisplaybreaks
\usepackage{graphicx}
\usepackage{xcolor}
\usepackage{comment}
\usepackage{array}
\usepackage{tikz-cd}
\usepackage{cite}
\usepackage{enumerate,enumitem,hyperref,cleveref}

\hypersetup{colorlinks,linkcolor={blue},citecolor={blue},urlcolor={blue}}

\newtheorem{theorem}{Theorem}[section]
\newtheorem{proposition}[theorem]{Proposition}
\newtheorem{lemma}[theorem]{Lemma}
\newtheorem{corollary}[theorem]{Corollary}
\theoremstyle{definition}

\theoremstyle{remark}
\newtheorem{remark}[theorem]{Remark}

\numberwithin{equation}{section}

\usepackage{color}
\newcommand{\R}{\mathbb R}
\newcommand{\dd}{\,d\mu_v}
\newcommand{\norm}[1]{\left\|#1\right\|}
\newcommand{\ip}[2]{\langle #1,#2\rangle}
\newcommand{\E}{\mathcal E}
\newcommand{\HH}{\mathcal H}
\newcommand{\one}{\mathbf1}
\DeclareMathOperator{\Dom}{Dom}

\title[Dimension-free weak $(1,1)$ estimates ...]
{Dimension-free weak $(1,1)$ estimates for Riesz transforms with constant drift}

\author[S. Mukherjee]{Suman Mukherjee}

\address{Department of Mathematics, Indian Institute of Technology Bombay, Powai, Mumbai--400076, India.}
\email{sumanmukherjee822@gmail.com}

\keywords{Riesz transform; Laplacian with drift; dimension-free estimate; partial balayage.}
\subjclass[2020]{Primary: 42B15, 42B20. Secondary: 35R11, 47G10}

\begin{document}

\begin{abstract}
Let $L_v=-\Delta-2v\cdot\nabla$ on $L^2(\mathbb R^n,e^{2v\cdot x}\,dx)$, where $v\in\mathbb R^n\setminus\{0\}$ is a constant drift vector. We prove that the full vector Riesz transform $\nabla L_v^{-1/2}$ satisfies a weak-type $(1,1)$ estimate for real-valued functions with constant at most $2$, uniformly in both the dimension and the drift vector. Our proof extends the fractional obstacle method of Ouyang, Spector and Stockdale (\href{https://arxiv.org/abs/2608.18068}{https://arxiv.org/abs/2608.18068}) to the exponentially weighted setting. Interpolation then yields strong $L^p$ bounds for $1<p\le 2$, with constants depending only on $p$.
\end{abstract}
\maketitle


\section{Introduction}\label{sec:intro}
Riesz transforms are among the basic examples of singular integral operators and provide a natural connection between first-order differentiation and the functional calculus of second-order operators. In Euclidean space, the vector transform
\[
 Rf:=\nabla(-\Delta)^{-1/2}f
   =\bigl(\partial_1(-\Delta)^{-1/2}f,\ldots,
          \partial_n(-\Delta)^{-1/2}f\bigr)
\]
is bounded on $L^p(\R^n)$ for $1<p<\infty$ and is of weak type $(1,1)$. These results belong to the classical theory of singular integrals developed by Calder\'on and Zygmund \cite{CZ}. The dependence of the operator norms on the dimension, however, is a separate question. A proof of boundedness in each fixed dimension need not provide an estimate uniform in $n$, and the endpoint is particularly sensitive to this distinction.

\subsection{Dimension-free estimates and the Euclidean endpoint}
The study of dimension-independent estimates for classical operators was developed systematically by Stein. For the Euclidean vector Riesz transform, Stein's work \cite{Stein83} gives
\begin{equation}\label{eq:classicalstrong}
 \norm{Rf}_{L^p(\R^n)}\le C_p\norm{f}_{L^p(\R^n)},
 \qquad 1<p<\infty,
\end{equation}
where $C_p$ is independent of $n$. Subsequent developments connected these inequalities with martingale methods and Bellman functions. In particular, Ba\~nuelos and Wang \cite{BW} obtained quantitative estimates through martingale inequalities, while Dragi\v{c}evi\'c and Volberg \cite{DV} developed Bellman-function proofs of dimensionless estimates for classical and Ornstein--Uhlenbeck Riesz transforms. These approaches demonstrated that the full vector can often be estimated without a loss from the number of its components.

The corresponding weak $(1,1)$ problem resisted the methods available for strong $L^p$ bounds. Stein formulated the dimension-free endpoint question in his 1986 International Congress address \cite[Problem b, p.~203]{Stein87}. For the individual Euclidean transforms, Janakiraman \cite{Jan} obtained a logarithmic dependence on the dimension, of order $\log (n)$ for $n\ge2$. Further work on the dimensional endpoint problem includes Spector and Stockdale \cite{SS}. In these questions, bounds for a single coordinate and bounds for the Euclidean length of the full vector must be distinguished: estimating the coordinates separately does not automatically give a dimension-free vector estimate.

Recently, Ouyang, Spector and Stockdale \cite{OSS} proved that, for real-valued $f$,
\begin{equation}\label{eq:OSS}
 \sup_{t>0}t\bigl|\{x\in\R^n:|Rf(x)|>t\}\bigr|
 \le2\norm{f}_{L^1(\R^n)}.
\end{equation}
Their proof introduces an obstacle decomposition adapted to the identity $R=\nabla(-\Delta)^{-1/2}$. For $f\ge0$, the decomposition has the form $f=m+(-\Delta)^{1/2}u$, where $m$ is bounded by a prescribed level, has the same mass as $f$, and equals that level on $\{u>0\}$. The transformed remainder is $\nabla u$, which vanishes outside this set. Thus the argument avoids estimating the Riesz kernel away from enlarged supports of bad functions, a step at which dimension dependence may enter a Calder\'on--Zygmund decomposition.

The variational structure has antecedents in obstacle problems and partial balayage. Classical variational inequalities are treated in Kinderlehrer and Stampacchia \cite{KS}; partial balayage on Riemannian manifolds is developed by Gustafsson and Roos \cite{GR}. The divisible-sandpile viewpoint provides another related formulation \cite{LP}. For nonlocal operators, the Lewy--Stampacchia estimates of Servadei and Valdinoci \cite{SV} give a relevant precedent for controlling the obstacle multiplier. The role of these ideas in \cite{OSS} is to produce a decomposition tailored to the Riesz transform rather than to a general singular kernel.

\subsection{Riesz transforms with drift}
An important direction in the theory replaces the Euclidean Laplacian by a symmetric diffusion operator on a weighted space. Bakry's work \cite{Bakry} connects Riesz-transform estimates with lower curvature bounds and provides a framework for dimension-free strong estimates. The constant-drift model considered here has vanishing Bakry--\'Emery curvature, but its natural measure has exponential growth. It therefore combines a simple differential expression with a geometry different from that of Lebesgue measure.

Let $v\in\R^n\setminus\{0\}$ be constant. The Laplacian with drift $v$ is given by
\begin{equation}\label{eq:objects}
 L_v=-\Delta-2v\cdot\nabla.
\end{equation}
Let $\mu_v$ be the measure given by $d\mu_v(x)=e^{2v\cdot x}\,dx$. Weighted integration by parts makes $L_v$ symmetric and nonnegative on $L^2(\mu_v)$. Riesz transforms for Laplacians with drift is then defined by
\[R_v=\nabla L_v^{-1/2}=\bigl(\partial_1L_v^{-1/2},\ldots,\partial_nL_v^{-1/2}\bigr).\]
Now, as the operator $L_v$ extends to a self-adjoint operator on
$L^2(\mu_v)$, the Riesz transform
$R_v$ is bounded on $L^2(\mu_v)$. Moreover, the Bakry--\'Emery curvature
on the manifold $(\mathbb{R}^{n}, |\cdot|, \mu_v)$ is known to vanish. This curvature
property yields dimension-independent $L^p$ bounds, for every $p>1$, for this
operator; see \cite[Theorem 4.1]{Bakry}. We also mention that $L^p$-boundedness results for this Riesz transform, together with extensions to more general classes of Riesz-type operators, were established
by Lohou\'e and Mustapha \cite{LM}. Li, Sj\"ogren and Wu \cite{LSW} proved weak type $(1,1)$ for the transform $R_v$, uniformly in $v$. Their statement does not give a dimension-independent constant. Our aim is therefore to establish dimension-free weak-type $(1,1)$ estimates for Riesz transforms with drift.

\subsection{Main results}
Our main result concerning the weak endpoint boundedness of $R_v$ is the following.
\begin{theorem}\label{thm:main}
For every $n\ge1$, every constant $v\in\R^n\setminus\{0\}$, and every real-valued $f\in L^1(\mu_v)$,
\begin{equation}\label{eq:weakmain}
 \norm{R_vf}_{L^{1,\infty}(\mu_v;\R^n)}\le2\norm{f}_{L^1(\mu_v)}.
\end{equation}
Here $R_v$ denotes the unique extension of its $L^2(\mu_v;\R^n)$ realization that is continuous from $L^1(\mu_v)$ to convergence in measure.
\end{theorem}
The estimate controls the entire gradient vector and is uniform both in dimension and in the magnitude and direction of the drift. It strengthens the quantitative conclusion of \cite[Theorem 1]{LSW} by giving a universal constant. The value $2$ is an upper bound obtained by the argument; its optimality is not asserted.

By the Marcinkiewicz interpolation theorem, together with the weak-type $(1,1)$ boundedness and the $L^2(\mu_v)$-isometry, we obtain the following strong-type boundedness result.
\begin{corollary}\label{cor:strong}
For $1<p<2$ and real-valued $f\in L^p(\mu_v)$,
\begin{equation}\label{eq:strongmain}
 \norm{R_vf}_{L^p(\mu_v;\R^n)}\le
 \left[4p\left(\frac1{p-1}+\frac1{2-p}\right)\right]^{1/p}\norm{f}_{L^p(\mu_v)}.
\end{equation}
At $p=2$, $\norm{R_vf}_{L^2(\mu_v;\R^n)}=\norm{f}_{L^2(\mu_v)}$. In particular, for every fixed $1<p\le2$, the strong operator norm is bounded independently of $n$ and $v$.
\end{corollary}
The strong boundedness in this corollary is already known in the drift setting. Its purpose here is to exhibit an explicit consequence of the endpoint theorem and the exact vector $L^2(\mu_v)$ identity. The displayed interpolation constant is not intended to be optimal.

The central ingredient is the following weighted obstacle decomposition. It is the constant-drift counterpart of the decomposition in \cite{OSS}.
\begin{theorem}\label{thm:obstacle}
Let $v\ne0$ be constant, $0\le f\in L^1(\mu_v)\cap L^2(\mu_v)$, and $\lambda>0$. There exist $ u\in L^1(\mu_v)\cap\Dom(L_v^{1/2})$ and $ m\in L^1(\mu_v)\cap L^\infty(\mu_v)$ such that
\[
f=m+L_v^{1/2}u \quad \text{in } L^2(\mu_v), \quad\text{and}
\]
\[
\textup{(i)}\; u\ge0,\ 0\le m\le\lambda;\qquad
\textup{(ii)}\; \norm{m}_{L^1(\mu_v)}=\norm{f}_{L^1(\mu_v)};\qquad
\textup{(iii)}\; \norm{m}_{L^2(\mu_v)}^2\le\lambda\norm{f}_{L^1(\mu_v)}.
\]
For $\Omega=\{u>0\}$, one has
\[
\textup{(i)}\; m=\lambda \ \text{a.e. on }\Omega;\qquad
\textup{(ii)}\; \mu_v(\Omega)\le\lambda^{-1}\norm{f}_{L^1(\mu_v)};\qquad
\textup{(iii)}\; \nabla u=0 \ \text{a.e. on }\Omega^c.
\]
\end{theorem}
Theorem \ref{thm:obstacle} retains the obstacle structure of \cite{OSS}, but its construction requires additional care in the exponentially weighted setting. The measure $\mu_v$ is not doubling, and the Fourier scaling argument used in \cite{OSS} to establish mass conservation does not apply directly. Although the spectral gap provides coercivity, the obstacle estimate initially yields only a bounded, locally integrable residual. The main difficulty is to recover its global integrability and exact weighted mass without introducing dimensional losses. 
\subsection{Method of proof and organization}
Write $A=L_v^{1/2}$ and let $\E$ be the closed quadratic form of $A$. For nonnegative data, we minimize
\[
 J(w)=\frac12\E(w)-\int fw\dd+\lambda\int w\dd
\]
over $w\ge0$ in $\Dom(A^{1/2})\cap L^1(\mu_v)$. The spectral gap $L_v\ge|v|^2I$ supplies coercivity. The heat kernel represents $\E$ by a positive integral of squared differences, from which the required truncation inequalities follow. A Lewy--Stampacchia argument gives a bounded residual $m$, initially only locally integrable.

The passage from local boundedness to exact conservation of mass is the principal issue introduced by the exponential measure. We construct expanding smooth cutoffs $\chi_R\rightarrow1$ satisfying $\norm{A\chi_R}_\infty\rightarrow 0$.
Pairing this estimate with the integrable minimizer proves $\int m\dd=\int f\dd$. Consequently $m\in L^2(\mu_v)$, and the weak equation upgrades $u$ from the form domain $\Dom(A^{1/2})$ to the operator domain $\Dom(A)$. Complementarity yields $m=\lambda$ on $\Omega$, while Sobolev locality gives $\nabla u=0$ on $\Omega^c$. These are the ingredients in Theorem~\ref{thm:obstacle}.

For signed data, we apply the decomposition to the positive and negative parts. The resulting exceptional set has measure at most $\norm{f}_{L^1(\mu_v)}/\lambda$, and the combined good part has squared $L^2$ norm at most $\lambda\norm{f}_{L^1(\mu_v)}$. Since $R_v(Au)=\nabla u$ and $\norm{R_vg}_{L^2(\mu_v)}=\norm{g}_{L^2(\mu_v)}$, choosing the obstacle level equal to the distribution threshold gives Theorem~\ref{thm:main}.

The remainder of the paper is organized as follows. Section~\ref{sec:setup} develops the weighted framework and collects the basic properties of the associated heat semigroup. Section~\ref{sec:fractional} introduces the fractional form and establishes the cutoff estimate. Section~\ref{sec:minimization} is devoted to the proof of the obstacle decomposition, including the multiplier estimate, conservation of mass, and the resulting regularity properties. Finally, Section~\ref{sec:endpoint} contains the proof of the weak boundedness theorem. Throughout the paper, we fix $n\ge1$ and a constant $v\ne0$, and write $a=|v|$ and $\rho_v(x)=e^{2v\cdot x}$. Unless stated otherwise, functions in the variational arguments are real-valued. The measure $\mu_v$ is finite on compact sets and has the same null sets as Lebesgue measure.

\section{The weighted Laplacian and its heat semigroup}\label{sec:setup}

\subsection{Weighted integration by parts}
For $g,h\in C_c^\infty(\R^n)$, integration by parts in each coordinate gives
\begin{align*}
 \int(-\partial_j^2g)h\rho_v\,dx
 &=\int(\partial_jg)\partial_j(h\rho_v)\,dx\\
 &=\int(\partial_jg)(\partial_jh)\rho_v\,dx
   +2v_j\int(\partial_jg)h\rho_v\,dx.
\end{align*}
The second term is canceled by the contribution of $-2v_j\partial_jg$. Summing in $j$ yields
\begin{equation}\label{eq:ibp}
 \int(L_vg)h\dd=\int\nabla g\cdot\nabla h\dd.
\end{equation}
This proves symmetry and nonnegativity on compactly supported smooth functions. The following conjugation specifies the self-adjoint realization without ambiguity.

\subsection{The ground-state calculation}
Define $Ug=e^{v\cdot x}g$. Then
\[
 \norm{Ug}_{L^2(dx)}^2=\int|g(x)|^2e^{2v\cdot x}\,dx=\norm{g}_{L^2(\mu_v)}^2.
\]
Its inverse is $U^{-1}h=e^{-v\cdot x}h$, so $U$ is unitary from $L^2(\mu_v)$ onto $L^2(dx)$. For smooth $h$,
\begin{align*}
 L_v(e^{-v\cdot x}h)
 &=e^{-v\cdot x}(-\Delta+a^2)h.
\end{align*}
Let $H_a:=-\Delta+a^2$ with its usual domain $H^2(dx)$ in $L^2(dx)$. We define the self-adjoint realization of the operator $L_v$ by
\begin{equation}\label{eq:conj}
 L_v=U^{-1}H_aU,\ \text{with } \Dom(L_v)=U^{-1}H^2(dx).
\end{equation}
Unitary conjugation preserves self-adjointness and positivity. Since $H_a\ge a^2I$, we have $L_v\ge a^2I$. This spectral gap will be used later.

\subsection{Fractional powers and their domains}
We use the unitary Fourier transform with phase $e^{-ix\cdot\xi}$ and normalization $(2\pi)^{-n/2}$. Then the multiplier of $H_a$ is $|\xi|^2+a^2$. Thus for $s\ge0$,
\[
 U L_v^s U^{-1}=H_a^s,\ \text{with }\Dom(H_a^s)=\left\{h\in L^2(dx):
    (|\xi|^2+a^2)^s\widehat h(\xi)\in L^2(d\xi)\right\}.
\]
In particular, with $A=L_v^{1/2}$ and $\HH=\Dom(A^{1/2})=\Dom(L_v^{1/4})$, we have $\Dom(A)=U^{-1}H^1(dx)$ and $\HH=U^{-1}H^{1/2}(dx)$. The equivalence with the usual Sobolev spaces follows because $a>0$ is fixed and $|\xi|^2+a^2$ is comparable to $1+|\xi|^2$. For $g, h \in \HH$, define
$$\E(g,h)=\ip{A^{1/2}g}{A^{1/2}h}_{L^2(\mu_v)}\ \text{and } \E(g)=\E(g,g).$$
The spectral gap implies
\begin{equation}\label{eq:gap}
 \E(g)\ge a\norm{g}_{L^2(\mu_v)}^2,\ \text{for all } g\in\HH.
\end{equation}
The norm $\bigl(\norm{g}_{L^2(\mu_v)}^2+\E(g)\bigr)^{1/2}$ makes $\HH$ a Hilbert space. We distinguish carefully between $\HH$, the domain of the form of $A$, and $\Dom(A)$, the domain of $A$ itself. We will see that the obstacle minimization initially gives only $u\in\HH$. The stronger conclusion $u\in\Dom(A)$ is proved later, after mass conservation.

\subsection{The weighted Sobolev identity}
Define
\[
 W^{1,2}(\mu_v)=\{g\in L^2(\mu_v):\partial_jg\in L^2(\mu_v),\ 1\le j\le n\},
\]
where derivatives are distributional derivatives on $\R^n$. This definition is meaningful because weighted $L^2$ functions are locally Lebesgue integrable. For $h=Ug$, distributional differentiation gives
\[
 U(\nabla g)=\nabla h-vh.
\]
It follows in both directions that $g\in W^{1,2}(\mu_v)$ if and only if $h\in H^1(dx)$. Hence
\begin{equation}\label{eq:domainA}
 \Dom(A)=W^{1,2}(\mu_v).
\end{equation}
For $h\in C_c^\infty(\mathbb{R}^n)$,
\begin{align*}
 \norm{\nabla U^{-1}h}_{L^2(\mu_v)}^2
 &=\int|\nabla h-vh|^2\,dx\\
 &=\int\bigl(|\nabla h|^2+a^2|h|^2\bigr)\,dx
      -\int v\cdot\nabla(h^2)\,dx\\
 &=\int\bigl(|\nabla h|^2+a^2|h|^2\bigr)\,dx.
\end{align*}
The last integral of a derivative is zero by compact support. Approximation in $H^1$ extends this identity to all $h\in H^1$. For complex $h$, we can replace $h^2$ by $|h|^2$ and the cross term by its real part. Plancherel therefore proves
\begin{equation}\label{eq:gradientidentity}
 \norm{\nabla g}_{L^2(\mu_v)}^2=\norm{Ag}_{L^2(\mu_v)}^2,
 \ \text{for all } g\in\Dom(A).
\end{equation}
By polarization, the corresponding bilinear identity also holds. This is the energy identity for $L_v$; it is different from the lower-order form $\E$ used in the obstacle functional.

\subsection{An operator core}
We will need to extend an identity, initially verified on $C_c^\infty(\mathbb{R}^n)$, to all of $\Dom(A)$. For this, we need an operator core. A subset of $\Dom(A)$ is called an operator core for $A$ if it is dense in $\Dom(A)$ with respect to the graph norm $\norm{g}_{L^2(\mu_v)}+\norm{Ag}_{L^2(\mu_v)}$.
\begin{proposition}\label{core}
 $C_c^\infty(\R^n)$ is an operator core for the operator $A$.   
\end{proposition}
\begin{proof}
Write $B=H_a^{1/2}$. By Plancherel theorem,
\[
\|h\|_{L^2(dx)}^2+\|Bh\|_{L^2(dx)}^2
 =\int(1+a^2+|\xi|^2)|\widehat h(\xi)|^2\,d\xi.
\]
This is equivalent to the $H^1(dx)$ norm squared. The standard cutoff-and-mollification construction supplies $h_k\in C_c^\infty(\R^n)$ with $h_k\to h$ and $Bh_k\to Bh$ in $L^2(dx)$ for every $h\in\Dom(B)$. Set $g_k=U^{-1}h_k$. Multiplying a smooth compactly supported function by $e^{-v\cdot x}$ preserves smoothness and support, so $g_k\in C_c^\infty$. Unitarity gives $g_k\to U^{-1}h$ and $Ag_k=U^{-1}Bh_k\to U^{-1}Bh$ in  $L^2(\mu_v)$. This is precisely the required graph-norm density.
\end{proof}

\begin{remark}
In particular, test functions belong to $\Dom(A)\subset\HH$. They also belong to $\Dom(L_v)$, since $U\varphi$ is smooth compactly supported and therefore in $H^2$.
\end{remark}

\subsection{The heat semigroup}\label{sec:heat}
Let $P_t=e^{-tL_v}$. From \eqref{eq:conj},
\[
 P_t=U^{-1}e^{-tH_a}U=U^{-1}e^{-a^2t}e^{t\Delta}U.
\]
Writing $q_t(z)=(4\pi t)^{-n/2}e^{-|z|^2/(4t)}$, we obtain initially for bounded compactly supported $g$
\begin{align}
 P_tg(x)
 &=\int q_t(x-y)e^{-v\cdot x-a^2t+v\cdot y}g(y)\,dy\nonumber\\
 &=\int q_t(x+2tv-y)g(y)\,dy.\label{eq:heat}
\end{align}
Relative to the weighted measure, the same formula reads
\begin{equation}\label{eq:heatweighted}
 P_tg(x)=\int p_t(x,y)g(y)\,d\mu_v(y),\ \text{where }
 p_t(x,y):=q_t(x-y)e^{-v\cdot(x+y)-a^2t}.
\end{equation}
Three properties follow directly.

First, $p_t\ge0$, so $P_t$ preserves positivity. Second, $p_t(x,y)=p_t(y,x)$, so the kernel is symmetric relative to $\mu_v$, in agreement with self-adjointness. Third, the Gaussian in \eqref{eq:heat} integrates to one, giving
\begin{equation}\label{eq:rows}
 \int p_t(x,y)\,d\mu_v(y)=1.
\end{equation}
Symmetry then also gives 
\begin{equation}\label{eq:rows2}
 \int p_t(x,y)\,d\mu_v(x)=1.   
\end{equation}
Hence $P_t1=1$ in the bounded-function sense. For $g\ge0$, Tonelli yields
\begin{equation}\label{eq:invariance}
 \int P_tg\dd=\int g\dd,
\end{equation}
allowing either integral to be infinite. For $g\in L^1$, positivity and \eqref{eq:invariance} imply
\[
 \norm{P_tg}_{L^1(\mu_v)}\le\int P_t|g|\dd=\norm{g}_{L^1(\mu_v)}.
\]
For bounded $g$, \eqref{eq:rows} gives $|P_tg(x)|\le\norm{g}_\infty$, so $P_t$ is also an $L^\infty$ contraction. Jensen's inequality and invariance similarly give the $L^2$ contraction property.

\begin{remark}
There is no contradiction between $P_t1=1$ and the $L^2$ spectral gap $L_v\ge a^2I$. The constant function $1$ is not in $L^2(\mu_v)$. In fact, after a rotation the weight grows exponentially in one coordinate, and the total measure is infinite. Conservativity is a statement about the kernel and bounded functions, not the existence of a zero eigenvector in $L^2$.
\end{remark}

\section{Fractional forms and approximation of constants}\label{sec:fractional}
\subsection{Fractional-power representation}
For $s>0$, substitute $r=ts$ and integrate by parts to obtain
\begin{align*}
 \int_0^\infty(1-e^{-ts})\frac{dt}{t^{3/2}}
 &=s^{1/2}\int_0^\infty(1-e^{-r})r^{-3/2}\,dr\\
 &=2s^{1/2}\int_0^\infty e^{-r}r^{-1/2}\,dr
 =2\sqrt\pi\,s^{1/2}.
\end{align*}
The boundary term is zero, since $(1-e^{-r})r^{-1/2}$ tends to zero both at zero and at infinity. The formula also holds at $s=0$. Thus
\begin{equation}\label{eq:scalar}
 s^{1/2}=\frac1{2\sqrt\pi}\int_0^\infty(1-e^{-ts})\frac{dt}{t^{3/2}}.
\end{equation}
This fixes the normalization used in every subsequent form identity.

\subsection{The fractional Dirichlet form}
Let $E_{L_v}$ be the spectral resolution of $L_v$, and for $g\in L^2(\mu_v)$, set $d\nu_g(s)=d\ip{E_{L_v}(s)g}{g}_{L^2(\mu_v)}$. Since $\ip{(I-P_t)g}{g}_{L^2(\mu_v)}=\int(1-e^{-ts})\,d\nu_g(s)\ge0$, Tonelli and \eqref{eq:scalar} give
\begin{equation}\label{eq:formdiag}
 \frac1{2\sqrt\pi}\int_0^\infty\ip{(I-P_t)g}{g}_{L^2(\mu_v)}\frac{dt}{t^{3/2}}
 =\int s^{1/2}\,d\nu_g(s).
\end{equation}
The right side is finite exactly when $g\in\HH$, and then equals $\E(g)$. The assertion with possible infinite value for general $g\in L^2(\mu_v)$ is useful as it allows us to prove that a truncation belongs to the form domain, instead of assuming that fact.

For a real $g\in L^2(\mu_v)$, symmetry and the identities \eqref{eq:rows} and \eqref{eq:rows2} imply
\begin{align*}
 &\iint(g(x)-g(y))^2p_t(x,y)\,d\mu_v(x)d\mu_v(y)\\
&=2\norm{g}_{L^2(\mu_v)}^2-2\ip{P_tg}{g}_{L^2(\mu_v)}
 =2\ip{(I-P_t)g}{g}_{L^2(\mu_v)}.
\end{align*}
All the terms in LHS are integrable. In fact, the squared difference is bounded by $2|g(x)|^2+2|g(y)|^2$, and the marginal identities \eqref{eq:rows} and \eqref{eq:rows2} gives
$$\iint(g(x)-g(y))^2p_t(x,y)\,d\mu_v(x)d\mu_v(y)\le 4\|g\|^2_{L^2(\mu_v)}.$$
Combining the last equality with \eqref{eq:formdiag} proves
\begin{equation}\label{eq:jumpdiag}
 \E(g)=\frac1{4\sqrt\pi}\int_0^\infty\iint
 (g(x)-g(y))^2p_t(x,y)\,d\mu_v(x)d\mu_v(y)\frac{dt}{t^{3/2}}.
\end{equation}
Again, finite energy is equivalent to membership in $\HH$.

For $g,h\in\HH$, apply Cauchy--Schwarz with respect to the positive measure occurring in \eqref{eq:jumpdiag}. It shows that the product of the two differences is absolutely integrable. Polarizing \eqref{eq:jumpdiag} gives
\begin{align}
\E(g,h)&=\frac1{4\sqrt\pi}\int_0^\infty\iint
 (g(x)-g(y))(h(x)-h(y))p_t(x,y)
 \,d\mu_v(x)d\mu_v(y)\frac{dt}{t^{3/2}}.\label{eq:jump}
\end{align}
In particular, this bilinear integral is an ordinary absolutely convergent integral, not a principal value.

\subsection{Contractions, disjoint supports, and negative parts}

The next lemma establishes the Markov property of the energy. More precisely, composing a function with any $1$-Lipschitz contraction that fixes the origin does not increase its energy. It also shows that if two nonnegative functions have disjoint supports, then their interaction with respect to the fractional energy is nonpositive. As a consequence, the positive and negative parts of a function remain in the energy space and satisfy the corresponding sign inequalities.
\begin{lemma}\label{lem:trunc}
If $g\in\HH$ and $\Phi:\R\to\R$ is $1$-Lipschitz with $\Phi(0)=0$, then $\Phi(g)\in\HH$ and $\E(\Phi(g))\le\E(g)$. If $g,h\in\HH$ are nonnegative and $gh=0$ a.e., then $\E(g,h)\le0$. Consequently, for every real $z\in\HH$,
\begin{equation}\label{eq:negative}
 \E(z,z^-)\le-\E(z^-),\ \text{where } z^-=\max\{-z,0\}.
\end{equation}
\end{lemma}
\begin{proof}
Since $|\Phi(s)|\le|s|$, the composition is in $L^2(\mu_v)$. The pointwise inequality
$|\Phi(g(x))-\Phi(g(y))|\le|g(x)-g(y)|$, together with the domain criterion in \eqref{eq:jumpdiag}, proves membership in $\HH$ and the energy estimate.

For nonnegative functions with $gh=0$, the integrand in \eqref{eq:jump} simplifies a.e. to
\[
 (g(x)-g(y))(h(x)-h(y))=-g(x)h(y)-g(y)h(x)\le0.
\]
The absolute integrability already proved justifies integrating this sign inequality. Finally, $z^+$ and $z^-$ belong to $\HH$ by the first assertion and satisfy $z^+z^-=0$. Since $z=z^+-z^-$,
\[
 \E(z,z^-)=\E(z^+,z^-)-\E(z^-,z^-)\le-\E(z^-).
\]
\end{proof}
\begin{remark}
The sign in \eqref{eq:negative} is central. It is the source of the upper bound for the obstacle multiplier; a Hilbert-space energy without this positivity structure would not suffice.
\end{remark}
\subsection{Expanding cutoffs}\label{sec:cutoff}

We next establish a semigroup representation for the fractional operator $A$ on test functions. 
In addition to convergence in $L^2(\mu_v)$, the representation converges uniformly, which allows us to obtain $L^\infty$ estimates for $A\phi$.

\begin{proposition}
\label{prop:operatorfraction}
For every $\phi\in C_c^\infty(\mathbb{R}^n)$,
\begin{equation}\label{eq:operatorfraction}
 A\phi=\frac1{2\sqrt\pi}\int_0^\infty(\phi-P_t\phi)\frac{dt}{t^{3/2}}.
\end{equation}
The integral converges both in $L^2(\mu_v)$ and in the uniform norm. Moreover,
\begin{equation}\label{eq:uniformdifference}
 \norm{\phi-P_t\phi}_\infty
 \le\min\{2\norm{\phi}_\infty,\ t\norm{L_v\phi}_\infty\}.
\end{equation}
\end{proposition}

\begin{proof}
For $\phi\in C_c^\infty(\mathbb{R}^n)$, the scalar formula
\eqref{eq:scalar} and the spectral theorem give
\begin{equation*}
 A\phi=\frac1{2\sqrt\pi}\int_0^\infty(\phi-P_t\phi)
 \frac{dt}{t^{3/2}}.
\end{equation*}
Here the integral converges in $L^r(\mu_v)$, where $r=2, \infty$. To verify convergence directly,
for $0<t\le1$ we can use the semigroup identity
\begin{equation}\label{eq:duhamel}
 \phi-P_t\phi=\int_0^tP_sL_v\phi\,ds
\end{equation}
and the $L^r$ contraction to bound the norm by
$t\norm{L_v\phi}_{L^r(\mu_v)}$. For $t\ge1$, we can directly get the bound as
$2\norm{\phi}_{L^r(\mu_v)}$. These give integrability against
$t^{-3/2}dt$ at both endpoints. 
We now identify the $L^2(\mu_v)$-limit. For
$0<\varepsilon<M<\infty$, define
\[
T_{\varepsilon,M}\phi
=
\frac{1}{2\sqrt{\pi}}
\int_\varepsilon^M
(\phi-P_t\phi)\,\frac{dt}{t^{3/2}}.
\]
Thus the truncated operator $T_{\varepsilon,M}$ has multiplier
\[
F_{\varepsilon,M}(s)
=
\frac{1}{2\sqrt{\pi}}
\int_\varepsilon^M
(1-e^{-ts})\,\frac{dt}{t^{3/2}},
\qquad
s=|\xi|^2+a^2.
\]
By the scalar identity \eqref{eq:scalar} we have $F_{\varepsilon,M}(s)\rightarrow s^{1/2}$, as $\varepsilon\rightarrow0$ and $M\rightarrow\infty$. Moreover,
\[
0\le F_{\varepsilon,M}(s)\le s^{1/2}.
\]
Therefore, $\widehat{UT_{\varepsilon,M}\phi}(\xi)
=
F_{\varepsilon,M}(|\xi|^2+a^2)\,
\widehat{U\phi}(\xi)$ converges pointwise to $(|\xi|^2+a^2)^{1/2}\widehat{U\phi}(\xi)
=
\widehat{UA\phi}(\xi)$. Hence dominated convergence gives $T_{\varepsilon,M}\phi
\rightarrow
A\phi$ in $L^2(\mu_v)$.

It remains to identify the uniform limit with the same function.
Let $\ell$ denote the uniform limit of $T_{\varepsilon,M}\phi$.
On every compact set $K\subset\mathbb R^n$,
\[
\|T_{\varepsilon,M}\phi-\ell\|_{L^2(K,\mu_v)}
\le
\|T_{\varepsilon,M}\phi-\ell\|_\infty\,
\mu_v(K)^{1/2},
\]
so uniform convergence implies local $L^2(\mu_v)$-convergence to $\ell$. On the other hand, the same truncated integrals converge globally in
$L^2(\mu_v)$ to $A\phi$. By uniqueness of the local
$L^2(\mu_v)$-limit, $\ell=A\phi$ almost everywhere. This completes the proof.
\end{proof}

We next construct a family of smooth cutoff functions adapted to the operator $A$. 
The following lemma shows that these cutoffs approximate the constant function $1$ while their contribution to the energy vanishes as $R\to\infty$.
\begin{lemma}\label{lem:cutoff}
There is a family $\chi_R\in C_c^\infty(\mathbb{R}^n)$, $R\ge1$, with $0\le\chi_R\le1$ and $\chi_R\rightarrow 1$ pointwise, such that
\begin{equation}\label{eq:cutofflimit}
 \norm{A\chi_R}_\infty\longrightarrow0.
\end{equation}
Also for every $u\in\HH\cap L^1(\mu_v)$,
\begin{equation}\label{eq:boundaryzero}
 \E(u,\chi_R)\longrightarrow0.
\end{equation}
\end{lemma}
\begin{proof}
Choose a nonnegative smooth radial function $\chi$, nonincreasing as a function of the radius, equal to one on $B(0,1)$ and zero outside $B(0,2)$. Set $\chi_R(x)=\chi(x/R)$. Radial monotonicity gives $\chi_R\rightarrow1$. By differentiation,
\[
 L_v\chi_R(x)=-R^{-2}(\Delta\chi)(x/R)-2R^{-1}v\cdot(\nabla\chi)(x/R).
\]
Thus, writing $b_R=\norm{L_v\chi_R}_\infty$,
\begin{equation}\label{eq:bR}
 b_R\le R^{-2}\norm{\Delta\chi}_\infty
           +2aR^{-1}\norm{\nabla\chi}_\infty\longrightarrow0.
\end{equation}
Equations \eqref{eq:operatorfraction} and \eqref{eq:uniformdifference} show that, for every $T>0$,
\begin{align*}
 \norm{A\chi_R}_\infty
 &\le\frac1{2\sqrt\pi}\left[
 b_R\int_0^Tt^{-1/2}\,dt+2\int_T^\infty t^{-3/2}\,dt\right]\\
 &=\frac1{2\sqrt\pi}\left(2b_R\sqrt T+\frac4{\sqrt T}\right).
\end{align*}
If $b_R>0$, choose $T=2/b_R$ to obtain the explicit estimate
\begin{equation}\label{eq:cutoffexplicit}
 \norm{A\chi_R}_\infty\le\frac{2\sqrt2}{\sqrt\pi}\sqrt{b_R}.
\end{equation}
If $b_R=0$, let $T\to\infty$ in the preceding bound, obtaining $A\chi_R=0$. This proves \eqref{eq:cutofflimit} in all cases.

It remains to explain the form pairing. The spectral theorem implies that for $u\in\HH$ and $\phi\in\Dom(A)$,
\begin{equation}\label{eq:pairing}
 \E(u,\phi)=\ip{u}{A\phi}.
\end{equation}
One may verify this first for spectrally truncated $u$ and then pass to the limit in $\HH$. Applying it with $\phi=\chi_R$ and using the $L^1$ assumption on $u$, we have
\[
 |\E(u,\chi_R)|=\left|\int uA\chi_R\dd\right|
 \le\norm{u}_{L^1(\mu_v)}\norm{A\chi_R}_\infty\longrightarrow0.
\]
This proves \eqref{eq:boundaryzero}.
\end{proof}

\begin{remark}
The weighted volume of the support of $\chi_R$ can grow very rapidly. No estimate for that volume occurs above. The relevant bound is a uniform bound on $A\chi_R$, paired with the already finite number $\norm{u}_{L^1(\mu_v)}$. The constants in \eqref{eq:bR} may depend on $n$ and $v$, but they are used only to establish a limit equal to zero. They never multiply $\norm{f}_{L^1(\mu_v)}$ in the final weak estimate.
\end{remark}

\section{The weighted obstacle decomposition}\label{sec:minimization}
\subsection{Variational construction}
Fix $0\le f\in L^1(\mu_v)\cap L^2(\mu_v)$ and $\lambda>0$. In this section $f$ and $\lambda$ remain fixed. Define
\begin{equation}\label{eq:cone}
 K=\{w\in\HH\cap L^1(\mu_v):w\ge0\text{ a.e.}\}
\end{equation}
and the functional
\begin{equation}\label{eq:functional}
 J(w)=\tfrac12\E(w)-\ip{f}{w}_{L^2(\mu_v)}+\lambda\int w\dd,\qquad w\in K.
\end{equation}
The functional is finite on $K$. In fact $w\in\HH\subset L^2$ and $f\in L^2$ make the middle term finite by Cauchy--Schwarz, and $w\in L^1$ makes the last term finite. The constraint $w\ge0$ is the obstacle constraint. The term $\lambda\int w$ penalizes the size of a positive competitor; it will also ensure that the minimizer is integrable.

The next lemma establishes the coercivity of the functional $J$ on the set $K$.

\begin{lemma}\label{lem:coercive}
For every $w\in K$,
\begin{equation}\label{eq:coercive}
 J(w)\ge\tfrac14\E(w)+\lambda\norm{w}_{L^1(\mu_v)}-a^{-1}\norm{f}_{L^2(\mu_v)}^2.
\end{equation}
\end{lemma}
\begin{proof}
By Cauchy--Schwarz and \eqref{eq:gap},
\[
 \ip{f}{w}_{L^2(\mu_v)}\le\norm{f}_{L^2(\mu_v)}\norm{w}_{L^2(\mu_v)}
 \le a^{-1/2}\norm{f}_{L^2(\mu_v)}\E(w)^{1/2}.
\]
The elementary inequality $bc\le b^2/4+c^2$, with $b=\E(w)^{1/2}$ and $c=a^{-1/2}\norm{f}_{L^2(\mu_v)}$, gives
\[
 \ip{f}{w}_{L^2(\mu_v)}\le\tfrac14\E(w)+a^{-1}\norm{f}_{L^2(\mu_v)}^2.
\]
Inserting this into \eqref{eq:functional} and using $\int w\dd=\norm{w}_{L^1(\mu_v)}$, we get the required inequality.
\end{proof}
\begin{remark}
The dependence on $a^{-1}$ in this auxiliary coercivity estimate does not enter the endpoint constant. This estimate provides existence separately for each fixed $v\ne0$; the endpoint constant will use only exact mass conservation and the $L^2$ identity.
\end{remark}
We next prove the existence and uniqueness of the minimizer of $J$ over $K$.

\begin{proposition}\label{prop:minimizer}
There exists a unique $u\in K$ minimizing $J$ over $K$.
\end{proposition}
\begin{proof}
Because $0\in K$ and $J(0)=0$, the infimum $\gamma=\inf_KJ$ is at most zero. Equation \eqref{eq:coercive} bounds it below by $-a^{-1}\norm{f}_{L^2(\mu_v)}^2$, so $\gamma$ is finite. Choose a sequence $w_k\in K$ with $J(w_k)\to\gamma$. Dropping finitely many terms, suppose $J(w_k)\le1$. Then
\[
 \tfrac14\E(w_k)+\lambda\norm{w_k}_{L^1(\mu_v)}
 \le1+a^{-1}\norm{f}_{L^2(\mu_v)}^2.
\]
Thus both $\E(w_k)$ and $\norm{w_k}_{L^1(\mu_v)}$ are uniformly bounded. Then equation \eqref{eq:gap} also bounds $\norm{w_k}_{L^2(\mu_v)}$. Consequently $\{w_k\}$ is a bounded sequence in the Hilbert space $\HH$.

Now weak sequential compactness of bounded sets in a Hilbert space gives, after taking a subsequence, $w_k\rightharpoonup u$ in $\HH$. The inclusion $\HH\hookrightarrow L^2(\mu_v)$ is continuous, so the convergence also holds weakly in $L^2(\mu_v)$. For every nonnegative $\phi\in L^2(\mu_v)$,
\[
 \ip{u}{\phi}=\lim_k\ip{w_k}{\phi}\ge0.
\]
Taking $\phi=u^-\in L^2(\mu_v)$ gives $-\norm{u^-}_{L^2(\mu_v)}^2\ge0$, hence $u\ge0$ a.e.

For each ball $B(0,R)$, its indicator belongs to $L^2(\mu_v)$. Therefore
\[
 \int_{B(0,R)}u\dd
 =\lim_k\int_{B(0,R)}w_k\dd
 \le\liminf_k\norm{w_k}_{L^1(\mu_v)}.
\]
Letting $R\to\infty$ and using monotone convergence gives
\begin{equation}\label{eq:L1lsc}
 \norm{u}_{L^1(\mu_v)}\le\liminf_k\norm{w_k}_{L^1(\mu_v)}<\infty.
\end{equation}
Thus $u\in K$. In particular, no weak compactness assertion in $L^1(\dd)$ is required.

The map $w\mapsto A^{1/2}w$ is bounded from $\HH$ into $L^2(\mu_v)$. Weak convergence therefore implies weak convergence of $A^{1/2}w_k$ to $A^{1/2}u$ in $L^2(\mu_v)$. Lower semicontinuity of the Hilbert norm gives $\E(u)\le\liminf_k\E(w_k)$. Also $\ip{f}{w_k}_{L^2(\mu_v)}\to\ip{f}{u}_{L^2(\mu_v)}$ by weak $L^2(\mu_v)$ convergence. Combining these facts with \eqref{eq:L1lsc}, and the lower bound on each nonnegative term, yields
\[
 J(u)\le\liminf_kJ(w_k)=\gamma.
\]
Since $u\in K$, necessarily $J(u)=\gamma$.

For $w,z\in K$, the quadratic identity \eqref{eq:jumpdiag} and linearity of the other terms give
\[
 J\left(\frac{w+z}{2}\right)
 =\frac{J(w)+J(z)}2-\frac18\E(w-z).
\]
If $w,z$ both minimize, the left side is at least $\gamma$ and the first term on the right is $\gamma$. Hence $\E(w-z)=0$. Equation \eqref{eq:gap} implies $\norm{w-z}_{L^2(\mu_v)}=0$. Thus the minimizer is unique up to null sets.
\end{proof}

\subsection{The obstacle multiplier}\label{sec:LS}

We first derive the variational inequality satisfied by the minimizer $u$. 
The convexity of the admissible set $K$ also yields the corresponding energy identity. This energy identity will later prove that the multiplier vanishes
wherever $u$ is positive.

\begin{proposition}
\label{prop:VI}
The minimizer $u$ satisfies
\begin{equation}\label{eq:VI}
 \E(u,w-u)\ge\int(f-\lambda)(w-u)\dd,\ \text{for every }w\in K.
\end{equation}
Moreover,
\begin{equation}\label{eq:energy}
 \E(u)=\int(f-\lambda)u\dd.
\end{equation}
\end{proposition}

\begin{proof}
For $w\in K$ and $0\le t\le1$, convexity of $K$ gives
$u+t(w-u)\in K$. Expanding the energy and the linear terms exactly, we get
\begin{align*}
 J(u+t(w-u))-J(u)
 &=t\left[\E(u,w-u)-\int(f-\lambda)(w-u)\dd\right]
   +\tfrac{t^2}{2}\E(w-u).
\end{align*}
The left side is nonnegative. we divide by $t>0$ and let $t\rightarrow0$
to obtain
\begin{equation*}
 \E(u,w-u)\ge\int(f-\lambda)(w-u)\dd\ \text{for any }\in K.
\end{equation*}
The integral is absolutely convergent. In fact, we can use $L^2(\mu_v)$ norm for the term with $f$
and $L^1(\mu_v)$ norm for the constant term. Choosing $w=0$ gives
$\E(u)\le\int(f-\lambda)u\dd$; choosing $w=2u$ gives the reverse
inequality. Thus, the proof of \eqref{eq:energy} is complete.
\end{proof}

We next introduce the multiplier associated with the variational inequality.
The following proposition shows that it defines a positive distribution and hence a locally finite positive Radon measure.

\begin{proposition}
\label{prop:multiplier}
For real $\phi\in C_c^\infty(\mathbb{R}^n)$, define
\begin{equation}\label{eq:eta}
 \eta(\phi)=\E(u,\phi)-\int(f-\lambda)\phi\dd.
\end{equation}
Then $\eta$ is a positive distribution. Moreover, $\eta$ extends to a
locally finite positive Radon measure, also denoted by $\eta$.
\end{proposition}

\begin{proof}
We first verify that $\eta$ is a distribution. On a fixed compact set
containing the support of $\phi$, the weight is bounded above. Also, by
\eqref{eq:pairing} and \eqref{eq:gradientidentity},
\[
 \E(\phi)=\ip{\phi}{A\phi}\le\norm{\phi}_{L^2(\mu_v)}
 \norm{A\phi}_{L^2(\mu_v)}
 =\norm{\phi}_{L^2(\mu_v)}\norm{\nabla\phi}_{L^2(\mu_v)}.
\]
Thus $\E(\phi)^{1/2}$ is controlled by fixed-compact-set bounds on
$\phi$ and its first derivatives. The estimate
\[
 |\eta(\phi)|\le\E(u)^{1/2}\E(\phi)^{1/2}
              +\norm{f}_{L^2(\mu_v)}\norm{\phi}_{L^2(\mu_v)}
              +\lambda\norm{\phi}_{L^1(\mu_v)}
\]
proves continuity in the test-function topology.

If $\phi\ge0$, then $u+\phi\in K$. Inserting this competitor into
\eqref{eq:VI} gives $\eta(\phi)\ge0$. Therefore $\eta$ is a positive
distribution. To show why it is a locally finite measure, on a compact
set choose a smooth nonnegative $\psi$ equal to one on that set. For any
smooth $\phi$ supported there,
\[
 -\norm{\phi}_\infty\psi\le\phi\le\norm{\phi}_\infty\psi.
\]
The positivity of $\eta$ yields
\[|\eta(\phi)|\le\eta(\psi)\norm{\phi}_\infty.\]
This provides a continuous positive extension to compactly supported
continuous functions. The representation theorem for positive linear
functionals gives a locally finite positive Radon measure, also denoted
by $\eta$.
\end{proof}

We next prove a Lewy--Stampacchia type estimate.
\begin{proposition}\label{prop:LS}
The measure $\eta$ satisfies
\begin{equation}\label{eq:LSmeasure}
 0\le\eta\le(\lambda-f)^+\,d\mu_v.
\end{equation}
Consequently $\eta=h\,d\mu_v$ for a measurable function $h$ with
\begin{equation}\label{eq:hdensity}
 0\le h\le(\lambda-f)^+\le\lambda\quad\text{a.e.}
\end{equation}
\end{proposition}
\begin{proof}
Fix $\phi\in C_c^\infty(\mathbb{R}^n)$ with $\phi\ge0$ and $\varepsilon>0$. Set
\[
 z=u-\varepsilon\phi,\qquad w=z^+.
\]
We verify admissibility of $w$ before using it. Since $z\in\HH$, Lemma~\ref{lem:trunc} gives $z^+\in\HH$. Pointwise, $u\ge0$ and $\phi\ge0$ imply $0\le z^+\le u$, so $z^+\in L^1$ as well. Hence $w\in K$. The same lemma gives $z^-\in\HH$; moreover,
\begin{equation}\label{eq:zminus}
 0\le z^-=(\varepsilon\phi-u)^+\le\varepsilon\phi.
\end{equation}
In particular, $z^-$ is integrable and supported inside the support of $\phi$, up to null sets.

Since $z^+=z+z^-$, we have $w-u=-\varepsilon\phi+z^-$. Substitution into \eqref{eq:VI} gives
\[
 -\varepsilon\E(u,\phi)+\E(u,z^-)
 \ge-\varepsilon\int(f-\lambda)\phi\dd+\int(f-\lambda)z^-\dd.
\]
Rearranging, using the definition of $\eta$,
\begin{equation}\label{eq:LSstep}
 \varepsilon\eta(\phi)
 \le\E(u,z^-)+\int(\lambda-f)z^-\dd.
\end{equation}
The two terms on the right can be estimated separately. Since $u=z+\varepsilon\phi$, Lemma~\ref{lem:trunc} and Cauchy--Schwarz for the form give
\begin{align*}
 \E(u,z^-)
 &=\E(z,z^-)+\varepsilon\E(\phi,z^-)\\
 &\le-\E(z^-)+\varepsilon\E(\phi)^{1/2}\E(z^-)^{1/2}\\
 &\le\frac{\varepsilon^2}{4}\E(\phi).
\end{align*}
The last step is $-b^2+cb\le c^2/4$ with $b=\E(z^-)^{1/2}$ and $c=\varepsilon\E(\phi)^{1/2}$. For the other term, using \eqref{eq:zminus} and nonnegativity, 
\[
 \int(\lambda-f)z^-\dd
 \le\int(\lambda-f)^+z^-\dd
 \le\varepsilon\int(\lambda-f)^+\phi\dd.
\]
After substituting the above two inequalities in \eqref{eq:LSstep} and v by $\varepsilon$,
\[
 \eta(\phi)\le\frac\varepsilon4\E(\phi)+\int(\lambda-f)^+\phi\dd.
\]
Letting $\varepsilon\rightarrow0$ proves the upper bound on every nonnegative smooth test function. Positivity gives the lower bound. By approximation of compactly supported continuous nonnegative functions and regularity of Radon measures, this is the measure inequality \eqref{eq:LSmeasure}. The dominating measure has locally integrable density, and therefore the Radon--Nikodym theorem gives $\eta=h\,d\mu_v$ with $0\le h\le(\lambda-f)^+$. The final inequality follows from the assumption $f\ge0$.
\end{proof}

We set
\begin{equation}\label{eq:mdef}
 m=\lambda-h.
\end{equation}
Then \eqref{eq:hdensity} implies $0\le m\le\lambda$ and also $m\ge\min\{f,\lambda\}$. By \eqref{eq:eta}, for every $\phi\in C_c^\infty(\mathbb{R}^n)$,
\begin{equation}\label{eq:weakid}
 \int m\phi\dd=\int f\phi\dd-\E(u,\phi).
\end{equation}
At this stage $m$ is bounded and locally integrable. We have not yet proved $m\in L^1(\dd)$ or $m\in L^2(\dd)$, and \eqref{eq:weakid} is not yet being asserted as an equality $Au=f-m$ of $L^2(\dd)$ functions. The next section establishes these facts in the order needed.

\subsection{Mass conservation and the contact set}\label{sec:completion}

We first show that the mass of the function $m$ is exactly preserved. 
The cutoff functions from Lemma~\ref{lem:cutoff} allow us to pass to the limit in the weak identity and obtain the conservation law.

\begin{proposition}\label{prop:mass}
The function $m$ belongs to $L^1(\mu_v)$ and satisfies
\begin{equation}\label{eq:mass}
 \norm{m}_{L^1(\mu_v)}=\int m\dd=\int f\dd=\norm{f}_{L^1(\mu_v)}.
\end{equation}
\end{proposition}

\begin{proof}
Insert $\phi=\chi_R$ from Lemma~\ref{lem:cutoff} into \eqref{eq:weakid}. Each integral is finite since, the cutoff has compact support, $m$ is bounded, and $f$ is integrable. Thus
\begin{equation*}
 \int m\chi_R\dd=\int f\chi_R\dd-\E(u,\chi_R).
\end{equation*}
The minimization already gave $u\in L^1(\dd)\cap\HH$, so the last term tends to zero by \eqref{eq:boundaryzero}. Since $f,m\ge0$ and $\chi_R\rightarrow1$, monotone convergence applies to both integral terms, even though the integrability of $m$ is not yet known. The right side tends to the finite number $\int f\dd$. It follows simultaneously that $m$ is integrable and that
\begin{equation*}
 \norm{m}_{L^1(\mu_v)}=\int m\dd=\int f\dd=\norm{f}_{L^1(\mu_v)}.
\end{equation*}

\end{proof}

We next upgrade the weak formulation to a genuine operator identity. 
The conservation of mass obtained above first yields the required $L^2(\dd)$ integrability of $m$, after which the operator-core property and self-adjointness of $A$ give the desired regularity of $u$.

\begin{proposition}\label{prop:operatorid}
We have
\begin{equation}\label{eq:operatorid}
 u\in\Dom(A)=W^{1,2}(\mu_v)\ \text{and } Au=f-m.
\end{equation}
In particular, this gives the required regularity improvement from
$\Dom(L_v^{1/4})$ to $\Dom(L_v^{1/2})$.
\end{proposition}

\begin{proof}
Because $0\le m\le\lambda$,
\begin{equation}\label{eq:mL2}
 \norm{m}_{L^2(\mu_v)}^2=\int m^2\dd\le\lambda\int m\dd
 =\lambda\norm{f}_{L^1(\mu_v)}<\infty.
\end{equation}
Hence $g=f-m\in L^2(\dd)$ (and also $L^1(\dd)$). Rewriting \eqref{eq:weakid} using \eqref{eq:pairing} gives
\begin{equation}\label{eq:coreweak}
 \ip{u}{A\phi}=\ip{g}{\phi},
 \ \text{for any }\phi\in C_c^\infty(\mathbb{R}^n).
\end{equation}
For any $\psi\in\Dom(A)$, the operator-core property
(Proposition~\ref{core}) provides
$\phi_k\in C_c^\infty(\mathbb{R}^n)$ such that
$\phi_k\to\psi$ and $A\phi_k\to A\psi$ in $L^2(\dd)$. Both sides of
\eqref{eq:coreweak} therefore converge, giving
\[
 \ip{u}{A\psi}=\ip{g}{\psi},\ \text{for any }\psi\in\Dom(A).
\]
So the linear functional
$\psi\mapsto\ip{u}{A\psi}$ is bounded for the $L^2(\dd)$ norm, with bound
$\norm{g}_{L^2(\mu_v)}$. 
By the definition of the adjoint, this says $u\in\Dom(A^*)$ and
$A^*u=g$. Since $A$ is self-adjoint, we conclude
\[
 u\in\Dom(A)=W^{1,2}(\mu_v)\ \text{and } Au=f-m.
\]
\end{proof}

Finally, we derive the complementarity relation between $u$ and $m$ by comparing the spectral and variational expressions for the energy. 
This relation determines $m$ on the positivity set of $u$ and gives a bound for the measure of that set.

\begin{proposition}
\label{prop:complementarity}
The functions $u$ and $m$ satisfy the complementarity relation
\begin{equation}\label{eq:complementarity}
 \int(\lambda-m)u\dd=0.
\end{equation}
Moreover, if $\Omega=\{u>0\}$, then
\begin{equation}\label{eq:contact}
 m=\lambda\text{ a.e. on }\Omega,\qquad
 \lambda\mu_v(\Omega)=\int_\Omega m\dd
 \le\norm{f}_{L^1(\mu_v)}.
\end{equation}
\end{proposition}

\begin{proof}
We have $u\in\Dom(A)$, so the spectral identity gives
\[
 \E(u)=\ip{Au}{u}=\int(f-m)u\dd.
\]
On the other hand, the variational identity \eqref{eq:energy} gives
$\E(u)=\int(f-\lambda)u\dd$. Subtracting the two identities yields
\[
 \int(\lambda-m)u\dd=0.
\]
All integrals used here are absolutely convergent. In fact, $fu$ and $mu$ are
integrable by Cauchy--Schwarz, and $\lambda u$ is integrable because
$u\in L^1(\dd)$. The integrand in \eqref{eq:complementarity} is nonnegative.
Therefore it vanishes a.e.. Define the measurable set
$\Omega=\{u>0\}$ using any representative of $u$. Then \eqref{eq:contact} follows immediately.

\end{proof}
We also use the following Sobolev locality lemma.
\begin{lemma}\label{lem:locality}
If $w\in W^{1,2}_{\mathrm{loc}}(\R^n)$, then $\nabla w=0$ a.e. on $\{w=0\}$.
\end{lemma}
\begin{proof}
First we recall the chain rule for a $C^1$ function $\Phi$ with bounded derivative $\nabla(\Phi(w))=\Phi'(w)\nabla w$ locally in the distributional sense. To justify it, on a relatively compact ball $B$, approximate $w$ in $W^{1,2}_{\mathrm{loc}}(\R^n)$ by smooth functions $w_k$ obtained by local mollification, and take a subsequence converging pointwise. The Lipschitz bound gives $\Phi(w_k)\to\Phi(w)$ in $L^2(B, dx)$ on the ball. For gradients, we can split
\[
 \Phi'(w_k)\nabla w_k-\Phi'(w)\nabla w
 =\Phi'(w_k)(\nabla w_k-\nabla w)
  +(\Phi'(w_k)-\Phi'(w))\nabla w.
\]
The first term tends to zero in $L^2(B, dx)$ by boundedness of $\Phi'$, and the second does so by dominated convergence, continuity of $\Phi'$, and $\nabla w\in L^2(B, dx)$. Passing to the distributional derivative gives the chain rule.

Now take $\Phi_\varepsilon(s)=\varepsilon\arctan(s/\varepsilon)$. Then
\[
 |\Phi_\varepsilon(s)|\le\frac\pi2\varepsilon,\qquad
 \Phi_\varepsilon'(s)=\frac1{1+(s/\varepsilon)^2}.
\]
On each bounded ball $B$, $\Phi_\varepsilon(w)\to0$ in $L^2(B, dx)$, whereas dominated convergence and the chain rule give
\[
 \nabla\Phi_\varepsilon(w)
 =\frac{\nabla w}{1+(w/\varepsilon)^2}
 \longrightarrow\one_{\{w=0\}}\nabla w\quad\text{in }L^2(B, dx).
\]
Since the distributional derivative of the limiting zero function is zero, the latter limit must equal zero. Exhausting $\R^n$ by balls proves the claim.
\end{proof}
For the minimizer $u$, \eqref{eq:operatorid} gives $u,\nabla u\in L^2(\mu_v)$. On any compact set the weight $\rho_v(x)$ is bounded below by a positive constant, so $u\in W^{1,2}_{\mathrm{loc}}(dx)$. Applying Lemma~\ref{lem:locality} and using equivalence of null sets, we  conclude
\begin{equation}\label{eq:localu}
 \nabla u=0\quad\mu_v\text{-a.e. on }\Omega^c=\{u=0\}
 \quad\text{up to null sets}.
\end{equation}
\subsection{Proof of the obstacle decomposition} Having established all the necessary ingredients, we are now ready to prove the obstacle decomposition.
\begin{proof}[{\bf Proof of Theorem~\ref{thm:obstacle}}]
Proposition~\ref{prop:minimizer} constructs $u\in K$. Proposition~\ref{prop:LS} and \eqref{eq:mdef} construct $m$ with the stated pointwise bounds. Equations \eqref{eq:mass}, \eqref{eq:mL2}, and \eqref{eq:operatorid} give integrability, the $L^2(\dd)$ estimate, and the operator identity. Finally, \eqref{eq:contact} and \eqref{eq:localu} give the assertions on $\Omega$. 
\end{proof}

\section{Proof of the weak-type estimate}\label{sec:endpoint}
\begin{proof}[{\bf Proof of Theorem~\ref{thm:main}}]

The spectrum of $A$ lies in $[a,\infty)$, so $A^{-1}$ is bounded on $L^2(\dd)$ with norm at most $a^{-1}$. For $g\in L^2(\dd)$, functional calculus gives $A^{-1}g\in\Dom(A)$ and $AA^{-1}g=g$. Then we have,
\[
 R_vg=\nabla A^{-1}g\in L^2(\mu_v;\R^n).
\]
By \eqref{eq:gradientidentity},
\begin{equation}\label{eq:L2isometry}
 \norm{R_vg}_{L^2(\mu_v;\R^n)}^2
 =\norm{\nabla A^{-1}g}_{L^2(\mu_v;\R^n)}^2
 =\norm{AA^{-1}g}_{L^2(\mu_v)}^2
 =\norm{g}_{L^2(\mu_v)}^2.
\end{equation}
Thus the full vector transform has $L^2$ norm exactly one. Also, for $u\in\Dom(A)$,
\begin{equation}\label{eq:RA}
 R_v(Au)=\nabla A^{-1}Au=\nabla u.
\end{equation}
Both identities are $L^2(\dd)$ identities. In particular, applying $R_v$ to the decomposition theorem is justified by boundedness.

Let $f\in L^1(\dd)\cap L^2(\dd)$ be real. Define
\[
 f^+=\max\{f,0\},\qquad f^-=\max\{-f,0\}.
\]
Then $f=f^+-f^-$, both parts belong to $L^1\cap L^2$, and
\begin{equation}\label{eq:signedmass}
 \norm{f^+}_{L^1(\mu_v)}+\norm{f^-}_{L^1(\mu_v)}=\norm{f}_{L^1(\mu_v)}.
\end{equation}
Fix a threshold $t>0$. For flexibility, first perform both obstacle decompositions at a possibly different level $\lambda>0$, that is, 
\[
 f^+=m_++Au_+,\qquad f^-=m_-+Au_-.
\]
Let $\Omega_+=\{u_+>0\}$, $\Omega_-=\{u_->0\}$, and set
\begin{equation}\label{eq:signedpieces}
 m=m_+-m_-,\qquad u=u_+-u_-,\qquad\Omega=\Omega_+\cup\Omega_-.
\end{equation}
Theorem~\ref{thm:obstacle} and \eqref{eq:signedmass} give
\begin{equation}\label{eq:OmegaBound}
 \mu_v(\Omega)\le\mu_v(\Omega_+)+\mu_v(\Omega_-)
 \le\frac{\norm{f}_{L^1(\mu_v)}}{\lambda}.
\end{equation}
Since both $m_+$ and $m_-$ take values in $[0,\lambda]$, their difference takes values in $[-\lambda,\lambda]$. Hence $|m|\le\lambda$. In particular,
\begin{equation}\label{eq:goodL2}
 \norm{m}_{L^2(\mu_v)}^2\le\lambda\int|m|\dd
 \le\lambda\int(m_++m_-)\dd=\lambda\norm{f}_{L^1(\mu_v)}.
\end{equation}
Here we have used the equality of masses of $m_\pm$ and $f^\pm$ in the last equality. One must not replace $|m_+-m_-|$ by $m_++m_-$ as an equality: the good parts may overlap, even though $f^+$ and $f^-$ are disjoint.

Subtracting the two decompositions gives $f=m+Au$ in $L^2(\dd)$. Equations \eqref{eq:RA} and \eqref{eq:signedpieces} yield
\[
 R_vf=R_vm+\nabla u.
\]
Outside $\Omega$, both gradients $\nabla u_+$ and $\nabla u_-$ vanish a.e., so $\nabla u=0$ there. Therefore, up to a null set,
\begin{equation}\label{eq:inclusion}
 \{|R_vf|>t\}\subset\Omega\cup\{|R_vm|>t\}.
\end{equation}

Chebyshev's inequality, the isometry \eqref{eq:L2isometry}, and \eqref{eq:goodL2} give
\[
 \mu_v\{|R_vm|>t\}
 \le t^{-2}\norm{R_vm}_{L^2(\mu_v;\R^n)}^2
 =t^{-2}\norm{m}_{L^2(\mu_v)}^2
 \le\frac{\lambda}{t^2}\norm{f}_{L^1(\mu_v)}.
\]
Combining with \eqref{eq:OmegaBound} and \eqref{eq:inclusion},
\begin{equation}\label{eq:twoterms}
 \mu_v\{|R_vf|>t\}
 \le\left(\frac1\lambda+\frac{\lambda}{t^2}\right)\norm{f}_{L^1(\mu_v)}.
\end{equation}
The expression in parentheses is minimized at $\lambda=t$. Equivalently, writing $\lambda=ct$, it becomes $(c^{-1}+c)/t\ge2/t$, with equality at $c=1$. We are free to select the level $\lambda$, so \eqref{eq:twoterms} proves
\begin{equation}\label{eq:weakdense}
 t\mu_v\{|R_vf|>t\}\le2\norm{f}_{L^1(\mu_v)},
\end{equation}
for any $t>0$ and for any real$\ f\in L^1(\dd)\cap L^2(\dd)$. Taking the supremum in $t$ gives the weak bound on the initial dense class. 

Since $L^1(\mu_v)\cap L^2(\mu_v)$ is dense in $L^1(\mu_v)$, the weak-type estimate and the standard weak-type extension principle yield a unique extension $R_v:L^1(\mu_v)\rightarrow L^{1,\infty}(\mu_v;\R^n)$ satisfying
\[
\|R_v f\|_{L^{1,\infty}(\mu_v;\R^n)}
\le 2\|f\|_{L^1(\mu_v)}.
\]

\end{proof}


\subsection*{Declaration of generative AI and AI-assisted technologies in the manuscript preparation process}

During the preparation of this manuscript, the author used OpenAI's ChatGPT solely for language editing and to improve the readability and clarity of the mathematical exposition. All mathematical statements and proofs were independently verified by the author. The author reviewed and revised all AI-assisted content as necessary and takes full responsibility for the final version of the article.

\subsection*{Acknowledgments} 
The author is supported by Institute Postdoctoral Fellowship from IIT Bombay.

\subsection*{Data availability}
Data availability is not applicable.

\subsection*{Competing interests}
The author declares that he has no competing interests.

\end{document}